\documentclass[10pt,a4paper,oneside]{amsart}

\usepackage{amssymb,amscd}
\usepackage[latin1,utf8]{inputenc}
\usepackage{verbatim}

\def\Q{\mathbb{Q}}

\def\R{\mathbb{R}}
\def\C{\mathbb{C}}
\def\H{\mathbb{H}}

\def\a{\alpha}

\def\w{\omega}
\def\G{\Gamma}

\def\H{\mathbf{h}}
\def\HH{\mathbf{H}}
\def\A{\mathbf{a}}
\def\AA{\mathbf{A}}
\def\Y{\mathbf{y}}
\def\YY{\mathbf{Y}}
\def\G{\mathbf{g}}
\def\Q{\mathbf{Q}}
\def\R{\mathbf{R}}
\def\X{\mathbf{x}}
\def\XX{\mathbf{X}}

\def\ssum{\sum_{t=0}^{L-1}}
\def\sssum{\sum_{t=0}^{L-1}\sum_{t'=0}^{L'-1}}
\def\vec{\mathrm{vec}}

\newtheorem{prop}{Proposition}

\title{Filtre NAPES pour signaux \`a bruit}
\begin{document}
\begin{abstract}
We extend the APES ({\sl Amplitude and Phase Estimation}) method of spectral analysis to the case of non-quasi periodic signals like appears in OFDM.\end{abstract}

\maketitle

\begin{center}
{\sc Jean-Philippe PR\' EAUX}\footnote[1]{\noindent Research center of the French air force (CReA), F-13661 Salon de
Provence air, France}
\footnote[2]{Laboratoire d'Analyse Topologie et Probabilit\'es,  Universit\'e de Provence, 39 rue F.Joliot-Curie, F-13453
marseille cedex 13, France.
\\ \indent\ {\it E-mail }: \ preaux@cmi.univ-mrs.fr
\smallskip\\
{\it Mathematical subject classification.} 65T04\smallskip\\
Preprint version, 2009.}
\end{center}

\section*{Notations} Donn\'e un nombre complexe $z\in\C$ on note respectivement $z^*$ et $|z|$ le conjugu\'e et le
module de $z$. On note $i\in\C$ le nombre complexe v\'erifiant $i^2=-1$.
Pour une matrice (ou un vecteur) $M$ \`a coefficients complexes on note respectivement $M^T$, $M^*$ et $M^C$ ses matrices
transpos\'ee, transconjugu\'ee et conjugu\'ee complexe ($M^*=(M^C)^T$). Si $M$ a $P$ lignes et $Q$ colonnes, ses lignes
seront num\'erot\'ees de 0 \`a $P-1$ et ses colonnes de 0 \`a $Q-1$. On consid\`ere la norme euclidienne 
$\| M\|=(\sum_{j=0}^{P-1}\sum_{k=0}^{Q-1}|M_{jk}|^2)^{\frac{1}{2}}$.

\section{Le cas unidimensionnel}

Soient 3 entiers strictement positifs, $M,L,N$, avec $N=M+L-1$. On consid\`ere $\{y_0,y_1,\ldots,y_{N-1}\}$ une suite de
donn\'ees complexes \`a analyser, et pour 
$t=0,1,\ldots ,L-1$, on pose $\Y(t)=(y_t,\ldots,y_{t+M-1})^T$ $\in\C^{M}$. Soit
$\w\in [0,2\pi[$ la fr\'equence \'etudi\'ee ; pour tout entier $K>0$ on consid\`ere
$\A_K(\w)=(1,e^{i\w},\ldots,e^{i(K-1)\w})^T\in\C^K$.
L'\'ecriture s'av\`ere plus concise en introduisant les deux matrices $M\times L$ \`a coefficients complexes,
$$\YY=(\Y_0,\Y_1,\ldots ,\Y_{L-1})$$
$$\AA(\w)=\A_M(\w)\,\A_L^*(\w)$$
et l'on souhaite exprimer :
$$\YY=\a(\w)\AA(\w)+res$$
o\`u $\a(\w)\in\C$ d\'esigne l'amplitude de fr\'equence $\w$ du signal.

\subsection{Le filtre APES} Soit $\H(\w)\in \C^M$ le
vecteur contenant les coefficients du filtre en fr\'equence $\w$. Le filtre APES   (cf. \cite{gro}) est la solution
$\H(\w)$ au probl\`eme d'optimisation sous contrainte en $\H(\w)\in\C^M$ et $\a(\w)\in\C$ :
$$
\min_{\H(\w),\a(\w)}\ \|\,\H^*(\w)\,\YY -\a(\w)\A_L^T(\w)\|\qquad \text{avec} \H^*(\w)\,\A_M(\w)=1
$$
soit encore,
$$
\min_{\H(\w),\a(\w)}\ \sum_{t=0}^{L-1} |\H^*(\w)\,\Y(t)-\a(\w)e^{i\w t}|^2\qquad \text{avec } \H^*(\w)\,\A_M(\w)=1
$$
Le filtre $\H(w)$ est choisi de sorte que :\\
$(i)$ le signal trait\'e $\H^*(\w)\,\YY$ soit le plus proche possible au sens des moindres carr\'es de la sinuso\"\i de
discr\`ete $\a(\w)\A_L^T(\w)$ de fr\'equence $\w$ et d'amplitude
complexe $\a(\w)$,\\
$(ii)$ le filtre ne distort pas la sinuso\"\i de discr\`ete $\A_M(\w)$ de fr\'equence $\w$.\\

 Il admet pour solution :
$${\a}_{APES}(\w)=\H^*_{APES}(\w)\,\G(w)$$
$$\H_{APES}(\w)=\frac{{\Q}^{-1}(\w)\,\A_M(\w)}{\A_M^*(\w)\,{\Q}^{-1}(\w)\,\A_M(\w)}$$
avec :
$$\G(\w)=\frac{1}{L}\sum_{t=0}^{L-1}\Y(t)e^{-i\w t}=\frac{1}{L} \YY\,{\A_L^C(\w)}$$
$$\Q(\w)=\frac{1}{L}(\sum_{t=0}^{L-1}\Y(t)\,\Y^*(t)) - \G(\w).\G^*(\w)=\frac{1}{L} \YY\,\YY^* -\G(\w)\,\G^*(\w)$$

\subsection{APES g\'en\'eralis\'e aux signaux \`a bruit} On g\'en\'eralise cette m\'ethode au cas avec bruit. Consid\'erons une suite
$x_0,x_1,\ldots,x_{N-1}$ d'\'el\'ements dans $\C$, repr\'esentant le bruit. Pour $0< K \leq N$ on note
$\X_K=(x_0,x_1,\ldots,x_{K-1})^T$. On note $\XX_{ref}$ la matrice $M\times L$ dont l'\'el\'ement ligne $i$ colonne $j$ est
$x_{i+j}$. On souhaite trouver l'amplitude complexe v\'erifiant :
$$\YY=\a(\w)\,\XX_{ref}\circ \AA(\w)+res$$
o\`u $\circ$ d\'esigne le produit de Hadamard de matrices.
 On cherche le filtre qui v\'erifie les deux conditions suivantes :\\
-- $(i)$ Le signal filtr\'e $\H^*(\w)\YY$ soit le plus proche possible (au sens des moindres carr\'es) de la sinuso\"\i de
\`a bruit
$\a(\w)\,(\X_L\circ \A_L(\w))$,\\
-- $(ii)$ le filtre ne d\'eforme pas les sinuso\"\i des \`a bruit $\X_M\circ \A_M(\w)$.\\

Trouver un filtre $\H(\w)$ v\'erifiant les conditions $(i)$ et $(ii)$ nous ram\`ene au probl\`eme d'optimisation en $\a(\w)$
et $\H(\w)$ sous contrainte  suivant :
$$
\min_{\H(\w),\a(\w)}\ \|\H^*(\w)\,\YY-\a(\w)\,\X_L\circ \A_L(\w)\|\quad\text{avec}\quad\H^*(\w)\,(\X_M\circ\A_M(\w))=1\\
$$
soit encore :
$$\min_{\H(\w),\a(\w)}\ \sum_{t=0}^{L-1} |\H^*(\w)\,\Y(t)-\a(\w)x(t)e^{i\w t}|^2\quad\text{avec}\quad\
\H^*(\w)\,(\X_M\circ\A_M(\w))=1
\qquad (*)\\
$$

\begin{prop}[{filtre NAPES}]
Le probl\`eme d'optimisation $(*)$ ci-dessus a pour solution (lorsqu'elle existe) :
\begin{gather*}
{\a}_{NAPES}(\w)=\H^*_{NAPES}(\w)\,\G(w)\\
\H_{NAPES}(\w)=\frac{{\Q}^{-1}(\w)\,(\X_M\circ\A_M(\w))}{(\X_M\circ\A_M(\w))^*\,{\Q}^{-1}(\w)\,(\X_M\circ\A_M(\w))}
\end{gather*}
 avec :
 \begin{gather*}
 \G(\w)=\frac{1}{\|\X_L\|^2}\sum_{t=0}^{L-1}x^*(t)e^{-i\w t}\,\Y(t)=\frac{1}{\|\X_L\|^2}\YY\,{(\X_L\circ \A_L(\w))^C}\\
\Q(\w)=\frac{1}{\|\X_L\|^2}\left(\sum_{t=0}^{L-1}\Y(t)\,\Y^*(t)\right) -
\G(\w)\,\G^*(\w)=\frac{1}{\|\X_L\|^2}\YY\,\YY^* - \G(\w)\,\G^*(\w)
\end{gather*}
 avec $\|\X_L\|^2=\ssum
|x_t|^2$.
\end{prop}

\subsection{D\'emonstration}
Nous suivons dans ses grandes lignes l'argument figurant dans \cite{gro} concernant la m\'ethode APES en le g\'en\'eralisant
au cas avec bruit c'est \`a dire \`a la r\'esolution de $(*)$.
\\

\begin{proof} Pour plus de lisibilit\'e nous notons $\a=\a(\w)$, $\A=\A(\w)$, $\H=\H(\w)$,
$\G=\G(\w)$ et $\Q=\Q(\w)$. Arrangeons la quantit\'e $I(\a,\H)$ \`a minimiser :
$$
\begin{aligned}
I(\a,\H)&=\sum_{t=0}^{L-1}|\H^*\,\Y(t)-\a\, x(t)e^{i\w t}|^2\\
&=\sum_{t=0}^{L-1}\left(\H^*\,\Y(t)-\a\, x(t)e^{i\w t}\right) \left(\Y^*(t)\,\H-\a^*x^*(t)e^{-i\w t}  \right)\\
&=\H^*\,\ssum\Y(t)\,\Y^*(t)\,\H+\ssum|\a|^2|x(t)|^2\\
 &-\a^*\H\,\ssum x^*(t)e^{-i\w t}\Y(t)-\a\,\ssum x(t)e^{i\w
t}\Y^*(t)\,\H\\
\end{aligned}
$$
On pose $ \R=\frac{1}{\|\X_L\|^2}\ssum \Y(t)\,\Y^*(t) $ et $\G=\G(\w)$ comme donn\'e ci-dessus,
$$
\begin{aligned}
I(\a,\H)&= \|\X_L\|^2\,\H^*\,\R\,\H+\| \X_L\|^2\,\left( |\a|^2 -\a^*\H^*\,\G -\a\G^*\,\H\right)\\
&=\| \X_L\|^2\,\left(\H^*\,\R\,\H+ |\a -\H^*\,\G|^2-|\H^*\,\G|^2\right)\\
&=\|\X_L\|^2\,\left(\H^*\,(\R-\G\,\G^*)\,\H+ |\a -\H^*\,\G|^2\right)\\
\end{aligned}
$$
se minimise pour $\a=\H^*\,\G$. Il reste \`a d\'eterminer le $\H$ optimal. On s'est ramen\'e  au probl\`eme d'optimisation :
$$
\min_{\H} \H^*\Q\H\qquad\text{sous la contrainte } \H^*\,(\X_M\circ\A)=1
$$
avec $\Q=\R-\G\,\G^*$. Il s'agit d'un probl\`eme quadratique avec une contrainte \'egalitaire lin\'eaire, et il est facile de
v\'erifier en appliquant les multiplicateurs de Lagrange qu'il admet pour solution (lorsqu'elle existe) :
$$\H=\frac{\Q^{-1}\,(\X_M\circ\A)}{(\X_M\circ\A)^*\,\Q^{-1}(\X_M\circ\A)}$$
ce qui donne la solution optimale recherch\'ee.
\end{proof}

\section{Le cas bidimensionnel}

Soient 6 entiers strictement positifs  $M,L,N,M',L',N'$ avec $N=M+L-1$ et $N'=M'+L'-1$. Soit $\YY_{N,N'}$ une matrice
$N\times N'$ \`a coefficients dans $\C$ repr\'esentant une s\'erie de donn\'ees 2-dimensionnelles. Pour $t=0,1,\ldots,L-1$ et
$t'=0,1,\ldots,L'-1$ soit $\Y_{t,t'}$ le vectoris\'e de la sous-matrice $M\times M'$ de $\YY_{N,N'}$ d\'ebutant ligne $t$,
colonne $t'$ :
$$
\Y_{t,t'}=\vec\left(\begin{array}{ccc} \YY_{N,N'}(t,t')& \cdots & \YY_{N,N'}(t,t'+M'-1)\\
\vdots & & \vdots\\
\YY_{N,N'}(t+M-1,t') &\cdots &\YY_{N,N'}(t+M-1,t'+M'-1)
\end{array}
   \right)_{M\times M'}
$$
C'est un vecteur de $\C^{M+M'}$. Soit $\w\in [0,2\pi[$ et pour tout entiers $P,P'$,
$$
\A_{P,P'}(\w,\w')=\left(\begin{array}{llll}1 &e^{i\w'}&\cdots&e^{i(P'-1)\w'}\end{array}\right)^T \otimes
\left(\begin{array}{llll}1 &e^{i\w}&\cdots&e^{i(P-1)\w}\end{array}\right)^T
$$
o\`u $\otimes$ d\'esigne le produit de Kronecker ; $\A_{P,P'}(\w,\w')$ est un vecteur de $\C^{P\times P'}$.

Un filtre bidimensionnel est une matrice $M\times M'$ \`a coefficients complexes.

\subsection{Filtre APES bidimensionnel}

C'est le filtre $\H(\w,\w')$ qui avec l'amplitude complexe $\a(\w,\w')$ est solution du probl\`eme d'optimisation:
\begin{gather*}
\min_{\H(\w,\w'),\a(\w,\w')}\ \sssum \mid \vec^*(\H(\w,\w'))\,\Y_{t,t'}-\a(\w,\w')e^{i(\w t+\w't')}\mid ^2\\
\vec^*(\H(w,w'))\,\A_{M,M'}(\w,\w')=1
\end{gather*}

Il admet pour solution :
$$
\a_{APES}(\w,\w')=\vec^*(\H(\w,\w'))\,\G(\w,\w')
$$
$$
\vec(\H_{APES}(\w,\w'))=\frac{\Q^{-1}(\w,\w')\,\A_{M,M'}(\w,\w')}{\A_{M,M'}^*(\w,\w')\,\Q^{-1}(\w)\,\A_{M,M'}(\w,\w')}
$$
avec,
$$
\G(\w,\w')=\frac{1}{L\times L'}\sssum e^{-i\,(\w t+\w't')}\Y_{t,t'}
$$
$$
\Q(\w,\w')=\frac{1}{L\times L'}\sssum \Y_{t,t'}\,\Y^*_{t,t'}-\G(\w,\w')\,\G(\w,\w')^*
$$

\subsection{Filtre APES bidimentionnel pour signaux \`a bruit connu}
C'est le filtre $\H(\w,\w')$ qui avec $\a(\w,\w')$ est solution du probl\`eme d'optimisation :

$$
\min_{\H(\w,\w'),\a(\w,\w')}\ \sssum \mid \vec^*(\H)\,\Y_{t,t'}-\a(\w,\w')x(t,t')e^{i(\w t+\w't')}\mid ^2\qquad (**)$$
$$
\vec^*(\H(\w,\w'))\,(\X_{M,M'}\circ\A_{M,M'})=1
$$

avec $x(t,t')\in\C$ repr\'esentant le bruit en $t,t'$ et $\X_{M,M'}=(x(0,0),\ldots ,x(M-1,0),\ldots,x(0,M'-1),\ldots
,x(M-1,M'-1))^T$.

\begin{prop} Le probl\`eme d'optimisation $(**)$ pr\'ec\'edent
se minimise pour
$$\a_{NAPES}(\w,\w')=\vec(\H_{NAPES}(\w,\w'))^*\,\G(\w,\w)$$
$$\vec(\H_{NAPES}(\w,\w'))=\frac{\Q^{-1}(\w,\w')\,(\X_{M,M'}\circ\A_{M,M'}(\w,\w'))}
{(\X_{M,M'}\circ\A_{M,M'})^*\Q^{-1}(\w,\w')(\X_{M,M'}\circ\A_{M,M'})}$$

avec
$$\G(\w,\w')=\frac{1}{\| \X_{L,L'}\|^2}\sssum x^*(t,t')e^{-i(\w t+\w't')}\Y(t,t')$$
$$\Q(\w,\w')=\frac{1}{\|\X_{L,L'}\|^2}\sssum \Y(t,t')\,\Y^*(t,t') -\G(\w,\w')\,\G^*(\w,\w')$$
o\`u $\| \X_{L,L'}\|^2=\sssum | x(t,t')|^2$.
\end{prop}

\subsection{D\'emonstration} La preuve suit toujours le m\^eme argument.

\begin{proof} On note $\a=\a(\w,\w')$, $\H=\H(\w,\w')$ et $\A_{M,M'}=\A_{M,M'}(\w,\w')$.
$$
\begin{aligned}
&\phantom{=}\sssum |\vec(\H)^*\,\Y(t,t')-\a\, x(t,t')e^{i(\w t+\w't')}|^2\\
 &=\sssum\left(\vec(\H)^*\,\Y(t,t')-\a\,
x(t,t')e^{i(\w t+\w't')}\right) \times
\left(\Y^*(t,t')\,\vec(\H)-\a^*\,x^*(t,t')e^{-i(\w t+\w't')}  \right)\\
&=\vec(\H)^*\,\sssum\Y(t,t')\,\Y^*(t,t')\,\vec(\H)+\sssum|\a\,|^2|x(t,t')|^2\\
 &-\a^*\,\vec(\H)\,\sssum x^*(t,t')e^{-i(\w t+\w't')}\Y(t,t')-\a\,\sssum x(t,t')e^{i(\w
t+\w't')}\Y^*(t,t')\,\vec(\H)\\
\end{aligned}
$$

On pose $ \R=\frac{1}{\|\X_{L,L'}\|^2}\ssum \Y(t,t')\,\Y^*(t,t') $,\\
et $\G=\G(\w,\w')=\frac{1}{\| \X_{L,L'}\|^2}\sssum x^*(t,t')e^{-i(\w t+\w't')}\Y(t,t')$
$$
\begin{aligned}
&= \| \X_{L,L'}\|^2\,\left(\vec(\H)^*\,\R\,\vec(\H)+ |\a|^2 -\a^*\vec(\H)^*\G -\a\,
\G^*\,\vec(\H)\right)\\
&=\| \X_{L,L'}\|^2\,\left( \vec(\H)^*\,\R\,\vec(\H)+|\a -\vec(\H)^*\,\G|^2-|\vec(\H)^*\,\G|^2\right)\\
&=\| \X_{L,L'}\|^2\,\left(\vec(\H)^*\,\left(\R-\G\,\G^*\right)\,\vec(\H)+|\a
 -\vec(\H)^*\,\G|^2\right)\\
\end{aligned}
$$
se minimise pour $\a=\vec(\H)^*\,\G$. Il reste \`a d\'eterminer le $\H$ optimal. On s'est ramen\'e  au probl\`eme
d'optimisation :
$$
\min_{\H} \vec(\H)^*\Q\,\vec(\H)\qquad\text{sous la contrainte } \vec(\H)^*\,(\X_{M,M'}\circ\A_{M,M'})=1
$$
avec $Q=\R-\G\,\G^*$. On obtient :
$$\vec(\H)=\frac{\Q^{-1}\,(\X_{M,M'}\circ\A_{M,M'})}{(\X_{M,M'}\circ\A_{M,M'})^*\Q^{-1}
(\X_{M,M'}\circ\A_{M,M'})}$$ ce qui donne la solution optimale recherch\'ee.
\end{proof}

\section{NAPES appliqu\'e aux signaux parcellaires}
Soit $N$ un entier non nul ; consid\'erons une suite de donn\'ees $y_0,y_1,\ldots , y_{N-1}$ dont certaines sont
manquantes, que l'on partitionne en ses sous-suites de donn\'ees connues de longueur $N_{2k-1}$ ($k=1,2,\ldots$)
entrecoup\'ees par ses sous-suites de donn\'ees manquantes de longueur $N_{2k}$ ($k=1,2,\ldots$), de sorte que les $N_1$
premi\`eres donn\'ees sont accessibles, les $N_2$ suivantes sont manquantes, et ainsi de suite. Soit $P$ le nombre de
segments disjoints ; ainsi $\sum_{k=1}^P N_k=N$. Notons $\Y_a$ et $\Y_u$ les 2 vecteurs de dimensions respectives
$N_1+N_3+\cdots +N_P$ et $N_2+N_4+\cdots +N_{P-1}$ contenant respectivement la suite de donn\'ees connues et la suite de
donn\'ees inconnues dans $y_0,y_1,\ldots , y_{N-1}$. On consid\`ere aussi la suite des bruits $(x_n)_{0\leq n<N}=
x_0,x_1,\ldots ,x_{N-1}$. \\

On applique NAPES pour obtenir une estimation de $\a(\w)$ et $\H(\w)$ de la façon suivante : choisir une longueur de
filtre initiale $M_0$. Si on v\'erifie la condition suivante, qui nous permet de construire une matrice $\R$ de
covariance de $M_0$ lignes constitu\'ee des donn\'ees accessibles :
$$\sum_{k=1,3,\ldots,P} \max(0,N_k-M_0+1)>M_0$$
appliquer 1.a, sinon appliquer 1.b, ci-dessous.\smallskip\\
1.a. Soient $L_k=N_k-M_0+1$, et notons $J$ le sous-ensemble de $\{1,3,\ldots,P\}$ pour lequel $L_k>0$.
 On note, si $L_1>0$, $\X_1$ le vecteur dont les coordonn\'ees sont donn\'ees par la sous-suite
$x_0,x_1,\ldots,x_{L-1}$ de $(x_n)_n$ d\'ebutant au rang $0$ et de longueur $L_1$ (si $L_1>0$), $etc \ldots,$ et pour
tout $P\in J$, $\X_P$ sa sous-suite d\'ebutant au rang $\sum_{k\in J}^{k<P}L_k$ et de longueur $L_P$.

Appliquer NAPES aux segments de donn\'ees correspondant \`a $J$ \`a l'aide des red\'efinitions suivantes :
$$\YY(\w)=\frac{1}{\sum_{k\in J} \| \X_k\|^2}\sum_{k\in J}
\sum_{t=N_1+N_2+\cdots+N_{k-1}}^{N_1+N_2+\cdots+N_{k-1}+L_k}\Y(t)x(t)e^{-i\w t}\, ,$$
$$\R=\frac{1}{\sum_{k\in J} \| \X_k\|^2}\sum_{k\in J}
\sum_{t=N_1+N_2+\cdots+N_{k-1}}^{N_1+N_2+\cdots+N_{k-1}+L_k}\Y(t)\Y^*(t)$$

\noindent 1.b. Appliquer NAPES avec $M_0=N/2$ \`a la suite de donn\'ees $y_0,y_1,\ldots,y_{N-1}$ en posant que $\Y_u$ est
une suite de z\'eros.\\

\noindent 2. On d\'etermine ensuite une estimation de $\Y_u$ bas\'ee sur l'estimation initiale des spectres $\a(\w_k)$ et
$\H(\w_k)$ d\'etermin\'ee pr\'ec\'edemment.  Elle est obtenue \`a l'aide de la minimisation suivante :

$$\min_{\Y_u}\sum_{k=0}^{K-1}\sum_{t=0}^{L-1} |\H^*(\w_k)\Y(t)-\a(\w_k)x_te^{i\w_kt}|^2$$
Consid\'erons la matrice $L\times N$ \`a coefficients dans $\C$
: $$\HH_k=\left(\begin{array}{llll} \H^*(\w_k)&0&\cdots&0\\
0&\H^*(\w_k) & &0\\
0& &\ddots&0\\
0&\cdots&0 &\H^*(\w_k)
\end{array}\right)$$
et
$$z_k=\a(\w_k)\,\X_*\circ \left(\begin{array}{c}
1\\
e^{i\w_k}\\
\vdots\\
e^{i(L-1)\w_k}
\end{array}\right)\in \C^L$$

Le crit\`ere de minimisation devient :
$$\sum_{k=0}^{K-1}\| \H_k\left(\begin{array}{c}
y_0\\
y_1\\
\vdots\\
y_{N-1}\end{array}\right)-z_k\|^2$$

Consid\'erons les matrices $A_k$ et $B_k$ telles que
$$\H_k\left(\begin{array}{c}
y_0\\
y_1\\
\vdots\\
y_{N-1}\end{array}\right)=A_k\Y_a+B_k\Y_u$$ et soit :
$$d_k=z_k-A_k\Y_a$$

De sorte que le crit\`ere \`a minimiser en $\Y_u$ devienne :
$$\sum_{k=0}^{K-1} \| B_k\Y_u -d_k\|^2$$
minimis\'e pour :
$$\Y_u=(\sum_{k=0}^{K-1} B_k^*B_k)^{-1}(\sum_{k=0}^{K-1}B_k^*d_k)$$

\noindent 3. Une fois qu'un estim\'e de $\Y_u$ a \'et\'e trouv\'e, on re-estime le spectre en appliquant NAPES cette fois-ci \`a
$\Y_a$ et $\Y_u$. Cel\`a consiste \`a minimiser en $\H(\w_k)$ et $\a(\w_k)$ le crit\`ere suivant :
$$\sum_{k=0}^{K-1}\sum_{t=0}^{L-1} |\H^*(\w_k)\Y(t)-\a(\w_k)x(t)e^{i\w_kt}|^2$$
Qui se scinde en fait en $K$ probl\`emes d'optimisations distincts pour d\'eterminer $\H(\w_k)$ et $\a(\w_k)$, pour $k$
variant de $0$ \`a $K-1$.\\

L'algorithme d'optimisation consiste alors en :\\
\noindent {\tt Soit $\delta>0$ la marge d'erreur impos\'ee\\
Appliquer l'\'etape initiale 1 pour estimer $\H(\w_k)$ et $\a(\w_k)$ pour $k=0,1,\ldots,K-1$\\
FAIRE :\\
\indent Estimer $\Y_u$ \`a l'aide de l'\'etape 2.\\
\indent Estimer $\H(\w_k)$ et $\a(\w_k)$ \`a l'aide de l'\'etape 3.\\
TANT QUE Les valeurs ont chang\'e de plus que $\delta$. }\\

Il s'agit en fait d'une optimisation cyclique du probl\`eme suivant :
$$\min_{\Y_u,\{\a(\w_k),\H(\w_k)\} } \sum_{k=0}^{K-1}\sum_{t=0}^{L-1} |\H^*(\w_k)\Y(t)-\a(\w_k)x(t)e^{i\w_kt}|^2$$

\end{document}